\documentclass[12pt]{article}
\usepackage{amsmath,amssymb,amsthm}
\usepackage{graphicx}
\usepackage{tikz}
\usetikzlibrary{patterns}
\usepackage{geometry}
\usepackage{hyperref}
\usepackage{caption}
\usepackage{float}
\usepackage{yhmath}
\usepackage{gensymb}
\usepackage[T1]{fontenc}
\usepackage{lmodern}
\usepackage{amsthm}
\theoremstyle{definition}
\newtheorem{definition}{Definition}[section]
\theoremstyle{plain}
\newtheorem{lemma}[definition]{Lemma}
\theoremstyle{remark}
\newtheorem{remark}[definition]{Remark}

\title{Dissecting Circles to Prove a Square: \\ A Novel Geometric Proof of the Pythagorean Theorem Using Circular Segments and Area Decomposition}
\author{Luca Nathanael Chang}

\date{June 20, 2025}

\begin{document}
\maketitle
\begin{abstract}
The Pythagorean Theorem has been proved in hundreds of ways, yet it inspires fresh insights through geometry and trigonometry. In this paper, we offer a new proof based on three circles that circumscribe the sides of a right triangle. Rather than invoke coordinate geometry, the argument relies purely on classical Euclidean constructions, trigonometric identities independent of the theorem itself, and a careful analysis of the areas of circular segments.

The key idea is to evaluate the area of the semicircle built on the hypotenuse in two distinct ways: directly and as a combination of areas formed by overlapping circular segments and triangles constructed on the legs of the triangle, as shown in Figure \ref{fig:Semicircle_Conglomerate}. Thales’ Theorem, inscribed angle theorem, basic trigonometric identities, and segment area formulas all play a role in a derivation that is both elementary and rigorous.

To the author’s knowledge, this specific approach, which combines circular symmetry, angle decomposition, and area comparison, has not appeared in the prior literature, including Loomis' comprehensive catalog \cite{loomis} and the extensive database at Cut-the-Knot \cite{bogomolny}. As such, it provides both a new perspective on an ancient theorem and an example of how classical tools can still yield original insights.
\end{abstract}


\section{Introduction}
Few results in mathematics are as universally recognized as the Pythagorean Theorem. From ancient Babylon to modern classrooms, the identity $a^2 + b^2 = c^2$ for the right triangles has served as both a cornerstone of geometry and a canvas for mathematical creativity. Over the centuries, hundreds of proofs have been devised using algebra, rearrangement, vectors, and even tiling. Still, new proofs continue to emerge, each revealing different facets of this deceptively simple equation.

This paper offers a novel geometric approach grounded in the areas of circular segments. Consider a right triangle with sides $a$, $b$, and $c$, and imagine placing a semicircle on each side, using that side as the diameter. The resulting configuration gives rise to several distinct circular segments, whose areas can be described using calculations of sector areas, as well as triangle areas, basic trigonometry, and classical theorems such as those of Thales.

The crux of the argument lies in calculating the area of the semicircle on the hypotenuse in two ways: once using standard formulas and once by decomposing it into overlapping circular segments and triangle regions derived from the side semicircles. These two paths to the same result converge—quite literally—on the Pythagorean theorem: $a^2+b^2=c^2$.

Unlike prior circle-based proofs that rely on coordinate geometry or power-of-a-point arguments, this proof remains purely geometric and trigonometric in spirit. The construction is direct and the reasoning is independent of the theorem itself.

Throughout this paper, we express angle measures in degrees (rather than radians) for pedagogical clarity and to align with geometric intuition, particularly when discussing arcs, sectors, and classical theorems such as Thales. All sector and segment area calculations use degree-based formulas. Readers who prefer radian measures may convert accordingly. Note that the use of square brackets [] is used to signify taking the area of the region within the brackets. Ie, [$\triangle$ABC] is equivalent to Area($\triangle ABC$).

\subsection*{Preliminary Definitions}

To ensure clarity, we define several key geometric terms used throughout this paper.

\begin{definition}[Circle]
A \emph{circle} is the set of all points in a plane that are equidistant from a fixed point, called the \emph{center}. The common distance is called the \emph{radius}.
\end{definition}

\begin{definition}[Chord]
A \emph{chord} is a line segment whose endpoints lie on the circumference of a circle.
\end{definition}

\begin{definition}[Arc]
An \emph{arc} is a continuous portion of the circumference of a circle, defined by two endpoints and the direction between them. An arc's measure can be measured in degrees, and its measure is the same as the central angle between the two radii that connect the center of the circle to the two defined endpoints. Minor Arc are denoted using 2 endpoints ($\wideparen{AB}$); Major Arcs use 3 points ($\wideparen{ABC}$). 
\end{definition}

\begin{definition}[Sector]
A \emph{sector} of a circle is the region enclosed by two radii and the arc connecting their endpoints. The area of a sector with central angle $\theta$ (in degrees) and radius $r$ is given by $\dfrac{\theta}{360^\circ} \pi r^2$.
\end{definition}

\begin{definition}[Circular Segment]
A \emph{circular segment} is the region bounded by a chord and the arc it subtends. It is obtained by subtracting the area of the triangle formed by the chord and the circle's center from the corresponding sector.
\end{definition}

\section{Construction} \label{sec:Construction}

Consider a right triangle $\triangle ABC$ with a right angle at $C$. Let the lengths of sides $\overline{AB}$, $\overline{BC}$, and $\overline{AC}$ be $c$, $a$, and $b$ respectively, where $a \leq b$ and $a, b, c \geq 0$.

We construct three circles as follows:
\begin{itemize}
    \item Circle $D$ with diameter $\overline{AB}$ (length $c$)
    \item Circle $E$ with diameter $\overline{AC}$ (length $b$)
    \item Circle $F$ with diameter $\overline{BC}$ (length $a$)
\end{itemize}

Points $D$, $E$, and $F$ represent the centers of these circles, which are the midpoints of the respective sides of $\triangle ABC$.

Let the intersection of Circles $E$ and $F$ be denoted as point $G$. By Thales' theorem\footnote{Thales' Theorem does not depend on the Pythagorean Theorem. (Elements, Book III, Proposition 31)  \cite{euclid1956elements}}, m$\angle CGB = 90^{\degree}$ and m$\angle CGA = 90^{\degree}$, which shows that point $G$ lies on line $\overline{AB}$ and $\overline{CG}$ is perpendicular to $\overline{AB}$. This diagram is shown below in Figure \ref{fig:init_diagram}.
\begin{figure}[H]
    \centering
    \begin{tikzpicture}[scale=1.5]
        \coordinate (C) at (0,0);
        \coordinate (B) at (3,0);
        \coordinate (A) at (0,4);
        \draw[thick] (A) -- (B) -- (C) -- cycle;       
        \draw (C) -- ++(0,0.2) -- ++(0.2,0) -- ++(0,-0.2);       
        \coordinate (D) at (1.5,2); 
        \coordinate (E) at (0,2); 
        \coordinate (F) at (1.5,0); 
        \draw[blue] (D) circle (2.5); 
        \draw[red] (E) circle (2); 
        \draw[green] (1.5,0) circle (1.5); 
        \coordinate (G) at (1.92,1.44);
        \draw[dashed] (G) -- (C);
        \draw[dashed] (E) -- (G);
        \draw[dashed] (1.5,0) -- (G);
        \node at (A) [above left] {A};
        \node at (B) [below right] {B};
        \node at (C) [below left] {C};
        \node at (D) [below] {D};
        \node at (E) [left] {E};
        \node at (F) [below] {F};
        \node at (G) [above right] {G};
        \node at (E) [above right] {$b$};
        \node at (F) [above left] {$a$};
        \node at (D) [above right] {$c$};
    \end{tikzpicture}
    \caption{Diagram For Proof Construction with Circumscribed Circles on Each Side}
    \label{fig:init_diagram}
\end{figure}
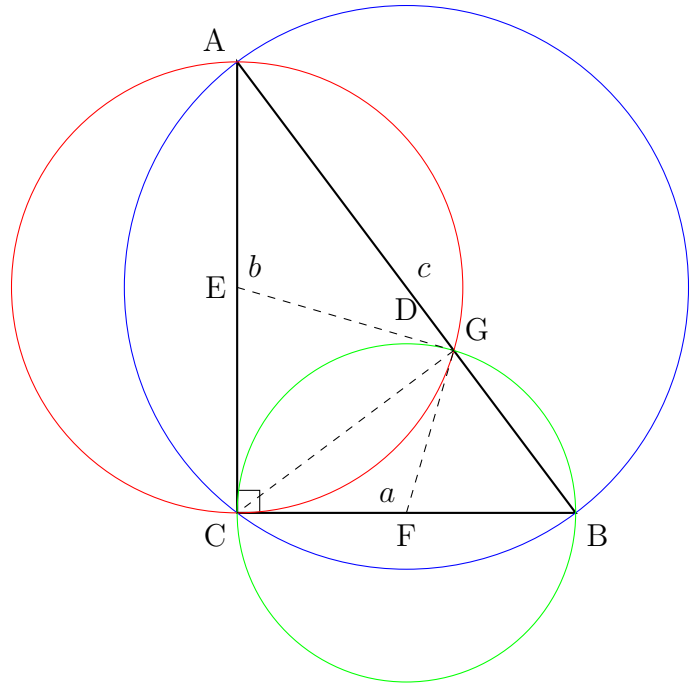

\section{Angle Relationships}
Let $m\angle CAB = \theta$. Utilizing Figure \ref{fig:init_diagram}, we can determine several angle relationships in our construction. Notably, since we have already established that $a \leq b$ and we know that the corresponding side to $\angle CAB$ has length $a$, we determine that the measure of $m\angle ABC \,\geq \,m\angle CAB$, as $b \geq a$. Since all triangles have angles that sum to $180^{\degree}$, and a right triangle has an angle with measure $90^{\degree}$, we determine that $m\angle CAB + m\angle ABC + 90^{\degree} = 180^{\degree} \:\Rightarrow\: m\angle CAB + m\angle ABC = 90^{\degree}$, which can be substituted to be $2 \cdot m\angle CAB \leq 90^{\degree} \:\Rightarrow\: m\angle CAB \leq 45^{\degree}\ \Rightarrow\: \theta\leq 45^{\degree}$. Other important relationships are as follows:

\begin{itemize}
    \item $\angle EAG = \angle EGA = \theta$, as $\overline{AE} = \overline{EG}$ (all radii of a circle are equal) and base angles of an isosceles triangle are congruent.
    \item $\angle CBG = 90^{\degree} - \theta$, as it is complementary to $\angle CAB$
    \item $\angle FBG = \angle FGB = 90^{\degree}-\theta$, as $\overline{FG} = \overline{FB}$ (all radii of a circle are equal) and base angles of an isosceles triangle are congruent.
    \item $\angle BFG = 180^{\degree} - \angle CBA - \angle FGB = 180^{\degree} - (90^{\degree}-\theta) - (90^{\degree}-\theta) = 2\theta$, as triangles have angles that sum to $180^{\degree}$.
    \item $\angle CFG = 180^{\degree} - \angle BFG = 180^{\degree}-2\theta$, as $\angle CFG$ and $\angle BFG$ are part of a linear pair, which sums to $180^{\degree}$.
    \item $\angle AEG = 180^{\degree} - \angle CAB - \angle EGA = 180^{\degree} - \theta - \theta = 180^{\degree} - 2\theta$, as a triangle has angles that sum to $180^{\degree}$.
    \item $\angle CEG = 180^{\degree} - \angle AEG = 180^{\degree} - (180^{\degree} - 2\theta) = 2\theta$,  as $\angle CEG$ and $\angle BEG$ are part of a linear pair, which sums to $180^{\degree}$.
\end{itemize}

\section{Line Segment Lengths}

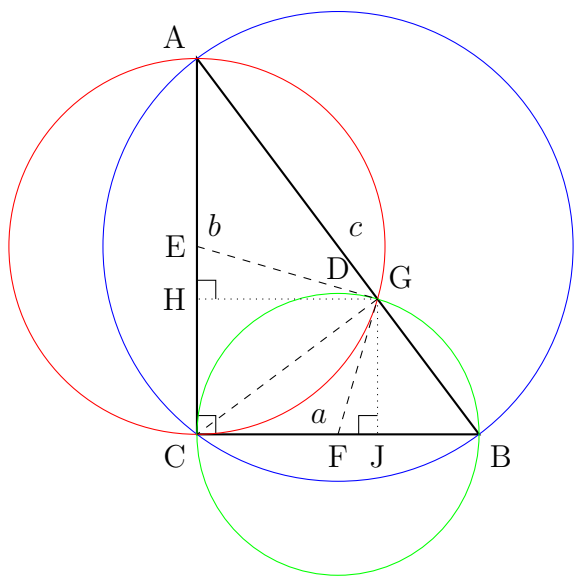
\begin{figure}[H]
    \centering
    \begin{tikzpicture}[scale=1.25]
        \coordinate (C) at (0,0);
        \coordinate (B) at (3,0);
        \coordinate (A) at (0,4);
        \coordinate (J) at (1.92, 0);
        \coordinate (H) at (0, 1.44);
        \draw[thick] (A) -- (B) -- (C) -- cycle;
        \draw (C) (0,0.2) -- ++(0.2,0) -- ++(0,-0.2);
        \draw (J) (1.72,0) -- ++(0,0.2) -- ++(0.2,0);
        \draw (H) (0,1.64) -- ++(0.2,0) -- ++(0,-0.2);
        \coordinate (D) at (1.5,2); 
        \coordinate (E) at (0,2); 
        \coordinate (F) at (1.5,0); 
        \draw[blue] (D) circle (2.5); 
        \draw[red] (E) circle (2); 
        \draw[green] (1.5,0) circle (1.5); 
        \coordinate (G) at (1.92,1.44);
        
        \draw[dashed] (G) -- (C);
        \draw[dashed] (E) -- (G);
        \draw[dashed] (1.5,0) -- (G);
        \draw[dotted] (G) -- (J);
        \draw[dotted] (G) -- (H);
        \node at (A) [above left] {A};
        \node at (B) [below right] {B};
        \node at (C) [below left] {C};
        \node at (D) [below] {D};
        \node at (E) [left] {E};
        \node at (F) [below] {F};
        \node at (G) [above right] {G};
        \node at (H) [left] {H};
        \node at (J) [below] {J};
        \node at (E) [above right] {$b$};
        \node at (F) [above left] {$a$};
        \node at (D) [above right] {$c$};
    \end{tikzpicture}
    \caption{Circumscribed Construction including Line Segments $\overline{GH}$ and $\overline{GJ}$}
    \label{fig:further_diagram}
\end{figure}

We draw altitudes from point $G$ to the sides $\overline{AC}$ and $\overline{CB}$ of $\triangle AGC$ and $\triangle CGB$ respectively, intersecting at points $H$ and $J$ respectively, as shown in Figure \ref{fig:further_diagram}. We note the sine double angle identity\footnote{While commonly derived using the unit circle, this can be proved using only geometric methods independent of the Pythagorean Theorem, making its use a valid, noncircular method to prove the Pythagorean theorem.\cite{bogomolny_sin}}:
$\sin(2\theta) = 2\cdot\sin(\theta)\cdot\cos(\theta)$. Noting that $\sin(\theta) = \dfrac{a}{c}$ and $\cos(\theta) = \dfrac{b}{c}$ from $\triangle ABC$ and the definition of sine and cosine as $\dfrac{opposite}{hypotenuse}$ and $\dfrac{adjacent}{hypotenuse}$ respectively, we determine:\footnote{Note that there is an alternate method that does not involve the use of trigonometry but keeps the same core idea: see \ref{sec:another method}}
\begin{align}
    \sin(2\theta) &= 2\cdot\sin(\theta)\cdot\cos(\theta) = 2 \cdot \dfrac{a}{c} \cdot \dfrac{b}{c} = \dfrac{2ab}{c^2}
\end{align}
From this, we can find $GJ$ and $GH$ by utilizing $\sin(\alpha)=\dfrac{opposite}{hypotenuse}$, which is rearranged to show that $opposite = \sin(\alpha) \cdot hypotenuse$:
\begin{align}
    GJ &= \sin(2\theta) \cdot FG =  \dfrac{2ab}{c^2} \cdot \dfrac{a}{2} = \dfrac{a^2b}{c^2} \\
    GH &= \sin(2\theta) \cdot EG = \dfrac{2ab}{c^2} \cdot \dfrac{b}{2} = \dfrac{ab^2}{c^2}
\end{align}

\section{Areas of Circular Segments}
\begin{figure}[H]
    \centering
    \begin{tikzpicture}[scale=1.1]
        \coordinate (C) at (0,0);
        \coordinate (B) at (3,0);
        \coordinate (A) at (0,4);
        \draw[thick] (A) -- (B) -- (C) -- cycle;
        \coordinate (D) at (1.5,2); 
        \begin{scope}
            \clip (D) circle (2.5); 
            \fill[blue!30, opacity=0.8, pattern=north east lines] (A) -- (C) arc (233.13:126.87:2.5) -- cycle;
        \end{scope}
        \coordinate (E) at (0,2); 
        \coordinate (F) at (1.5,0); 
        \draw[blue] (D) circle (2.5); 
        \draw[red] (E) circle (2); 
        \draw[green] (1.5,0) circle (1.5); 
        
        \draw[dashed] (C) -- (D);
        \draw[dashed] (A) -- (D);
        \node at (A) [above left] {A};
        \node at (B) [below right] {B};
        \node at (C) [below left] {C};
        \node at (D) [right] {D};
        \node at (E) [right] {$b$};
        \node at (F) [above] {$a$};
    \end{tikzpicture}
    \caption{Region A (shaded)}
    \label{fig:fig_a}
\end{figure}
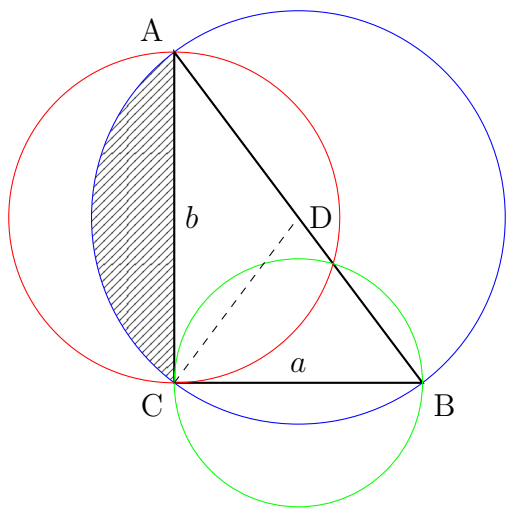

We calculate the areas of six key circular segments formed by our construction:
\subsection{\texorpdfstring{Region A: The area of the circular segment $\wideparen{AC}$ subtended by chord $\overline{AC}$ on circle $D$ (Figure \ref{fig:fig_a})}{Region A: The area of the circular segment subtended by chord AC on circle D (Figure 3)}}

By the inscribed angle theorem, the central angle $\angle CDA$ of Circle $D$ has an angle measure of $180^{\degree} - 2\theta$ since we have $\angle ABC = 90^{\degree} - \theta$. Thus, the measure of minor arc $\wideparen{CA} $ is $180^{\degree} - 2\theta$. This gives the area of sector $CDA$:
\begin{align}
    \text{Area of sector $[\wideparen{CDA}]$} &= \pi(\dfrac{c}{2})^2 \cdot\dfrac{(180^{\degree} - 2\theta)}{360^{\degree}} = \dfrac{\pi c^2}{4} \cdot \dfrac{180^{\degree} - 2\theta}{360^{\degree}}
\end{align}
The area of $[\triangle CDA]$ is $\dfrac{ab}{4}$, as it is half of the area of $[\triangle ABC]$ which has area $\dfrac{ab}{2}$. Thus, the area of the circular segment subtended by chord $\overline{AC}$ on circle $D$ is:
\begin{align}
    \text{Area of Circular Segment $[\wideparen{CA}]$} &= \dfrac{\pi c^2}{4} \cdot \dfrac{180^{\degree} - 2\theta}{360^{\degree}} - \dfrac{ab}{4} \hypertarget{eq:region a}{}
\end{align}

\subsection{\texorpdfstring{Region B: The area of the circular segment $\wideparen{CB}$ subtended by chord $\overline{CB}$ on circle $D$ (Figure \ref{fig:fig_b})}{Region B: The area of the circular segment subtended by chord CB on circle D (Figure 4)}}
\texorpdfstring{By the inscribed angle theorem,   central angle $\angle CDB$ has an angle of $2\theta$ due to $\angle CAB$ having measure $\theta$. Thus, the measure of minor arc $\wideparen{CB}$ of Circle $D$ is $2\theta$. From this, we find the area of sector $[\wideparen{CDB}]$:}{By the inscribed angle theorem, the central angle CDB has an angle of $2\theta$ due to CAB having measure $\theta$. Thus, the measure of minor arc CB of Circle D is $2\theta$. From this, we find the area of sector [CDB]:}
\begin{align}
    \text{Area of sector $[\wideparen{CDB}]$} &= \pi(\dfrac{c}{2})^2 \cdot \dfrac{2\theta}{360^{\degree}} = \dfrac{\pi c^2}{4} \cdot \dfrac{2\theta}{360^{\degree}}
\end{align}
Subtracting the area of $[\triangle CDB]$, $\dfrac{ab}{4}$, as it is $\dfrac{1}{2}$ of the area of $[\triangle ABC]$  yields:
\begin{align}
    \text{Area of Circular Segment $[\wideparen{CB}]$} &= \dfrac{\pi c^2}{4} \cdot \dfrac{2\theta}{360^{\degree}} - \dfrac{ab}{4} \hypertarget{eq:region b}{}
\end{align}
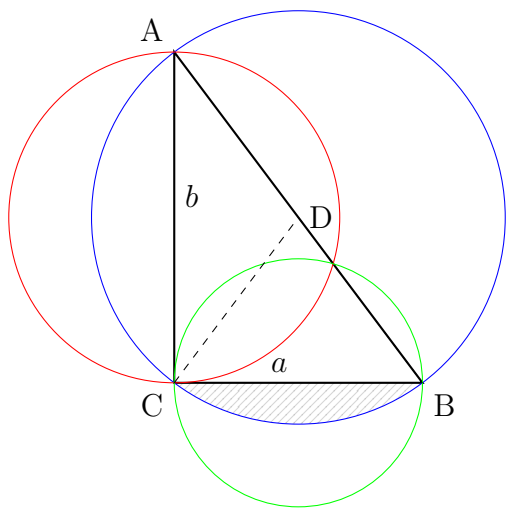
\begin{figure}[H]
    \centering
    \begin{tikzpicture}[scale=1.1]
        \coordinate (C) at (0,0);
        \coordinate (B) at (3,0);
        \coordinate (A) at (0,4);
        \draw[thick] (A) -- (B) -- (C) -- cycle;   
        \coordinate (D) at (1.5,2); 
        \begin{scope}
            \clip (D) circle (2.5); 
            \fill[blue!30, opacity=0.4, pattern=north east lines] (C) -- (B) arc (306.87:233.13:2.5) -- cycle;
        \end{scope}
        \coordinate (E) at (0,2); 
        \coordinate (F) at (1.5,0); 
        \draw[blue] (D) circle (2.5); 
        \draw[red] (E) circle (2); 
        \draw[green] (1.5,0) circle (1.5); 
        \draw[dashed] (C) -- (D);
        \draw[dashed] (D) -- (B);
        \node at (A) [above left] {A};
        \node at (B) [below right] {B};
        \node at (C) [below left] {C};
        \node at (D) [right] {D};
        \node at (E) [above right] {$b$};
        \node at (F) [above left] {$a$};  
    \end{tikzpicture}
    \caption{Region B (shaded)}
    \label{fig:fig_b}
\end{figure}

\subsection{\texorpdfstring{Region C: The area of the circular segment $\wideparen{AG}$ subtended by chord $\overline{AG}$ on circle $E$ (Figure \ref{fig:fig_c})}{Region C: The area of the circular segment subtended by chord AG on circle E (Figure 5)}}
We know that $m\angle AEG = 180^{\degree}-2\theta$. Thus, the area of the sector $[\wideparen{AEG}]$ is $\dfrac{\pi b^2}{4} \cdot \dfrac{180^{\degree}-2\theta}{360^{\degree}}$. The area of $[\triangle AEG]$ is $\dfrac{1}{2} GH \cdot AE$.  Substituting values, we obtain the area of $[\triangle AEG] = \dfrac{1}{2}(\dfrac{ab^2}{c^2} \cdot \dfrac{b}{2})$. Thus, we get the area below.
\begin{align}
    \text{Area of Region C: Upper-Right Circular Segment of Circle $E$}\nonumber 
\end{align}
\begin{align}
    =\dfrac{\pi b^2}{4} \cdot \dfrac{180^{\degree} - 2\theta}{360^{\degree}} - \dfrac{1}{2}(\dfrac{ab^2}{c^2} \cdot \dfrac{b}{2}) \hypertarget{eq:region c}{}
\end{align}
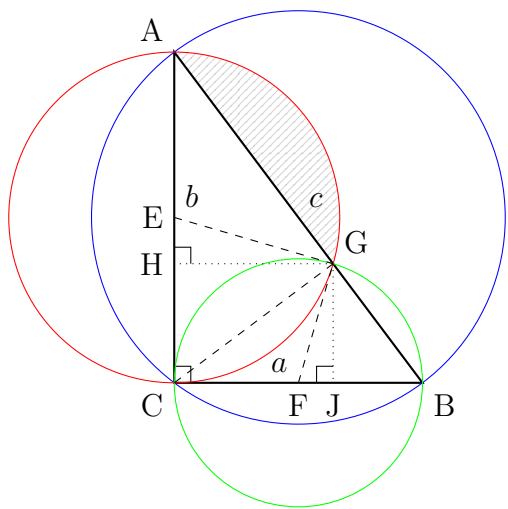
\begin{figure}[H]
    \centering
    \begin{tikzpicture}[scale=1.1]
        \coordinate (C) at (0,0);
        \coordinate (B) at (3,0);
        \coordinate (A) at (0,4);
        \draw[thick] (A) -- (B) -- (C) -- cycle;
        \draw (C) (0,0.2) -- ++(0.2,0) -- ++(0,-0.2);
        \draw (J) (1.72,0) -- ++(0,0.2) -- ++(0.2,0);
        \draw (H) (0,1.64) -- ++(0.2,0) -- ++(0,-0.2);
        \coordinate (D) at (1.5,2); 
        \coordinate (E) at (0,2); 
        \begin{scope}
            \clip (0,2) circle (2);
            \fill[red!30, opacity=0.4, pattern=north east lines]
                (0,4) -- (1.92,1.44)
                plot[domain=-16.26:90, samples=100]
                    ({0 + 2*cos(\x)}, {2 + 2*sin(\x)})
                -- cycle;
        \end{scope}
        \coordinate (F) at (1.5,0); 
        \draw[blue] (D) circle (2.5); 
        \draw[red] (E) circle (2); 
        \draw[green] (1.5,0) circle (1.5); 
        \coordinate (G) at (1.92,1.44);
        \coordinate (J) at (1.92, 0);
        \coordinate (H) at (0, 1.44);
        \draw[dashed] (G) -- (C);
        \draw[dashed] (E) -- (G);
        \draw[dashed] (1.5,0) -- (G);
        \draw[dotted] (G) -- (J);
        \draw[dotted] (G) -- (H);
        \node at (A) [above left] {A};
        \node at (B) [below right] {B};
        \node at (C) [below left] {C};
        \node at (E) [left] {E};
        \node at (F) [below] {F};
        \node at (G) [above right] {G};
        \node at (H) [left] {H};
        \node at (J) [below] {J};
        \node at (E) [above right] {$b$};
        \node at (F) [above left] {$a$};
        \node at (D) [above right] {$c$};
    \end{tikzpicture}
    \caption{Region C (shaded)}
    \label{fig:fig_c}
\end{figure}

\subsection{\texorpdfstring{Region D: The area of the circular segment $\wideparen{CG}$ subtended by chord $\overline{CG}$ on circle $E$ (Figure \ref{fig:fig_d})}{Region D: The area of the circular segment subtended by chord CG on circle E (Figure 6)}}
We know that $m\angle CEG = 2\theta$. Thus, the area of sector $[\wideparen{CEG}]$ is $\dfrac{\pi b^2}{4} \cdot \dfrac{2\theta}{360^{\degree}}$. The area of  $[\triangle CEG]$ is $\dfrac{1}{2} GH \cdot EC$.  Substituting values, we see that the area of $[\triangle CEG] = \dfrac{1}{2}(\dfrac{ab^2}{c^2} \cdot \dfrac{b}{2})$.  Thus, we get the area below.
\begin{align}
    \text{Area of Region D: Bottom-Right Circular Segment of Circle $E$}\nonumber
\end{align}
\begin{align}
    &= \dfrac{\pi b^2}{4} \cdot \dfrac{2\theta}{360^{\degree}} - \dfrac{1}{2}(\dfrac{ab^2}{c^2} \cdot \dfrac{b}{2})  \hypertarget{eq:region d}{}
\end{align}
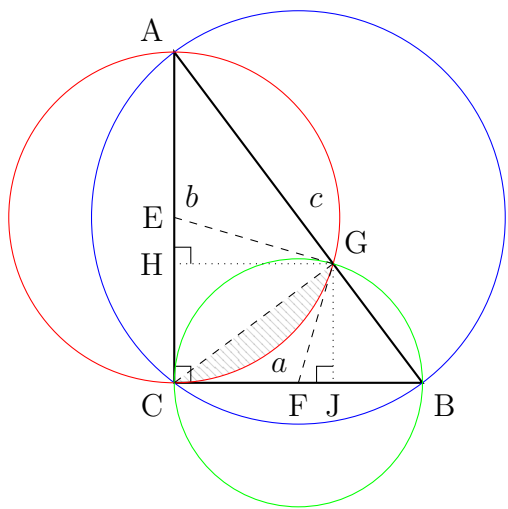
\begin{figure}[H]
    \centering
    \begin{tikzpicture}[scale=1.1]
        \coordinate (C) at (0,0);
        \coordinate (B) at (3,0);
        \coordinate (A) at (0,4);
        \draw[thick] (A) -- (B) -- (C) -- cycle;
        \draw (C) (0,0.2) -- ++(0.2,0) -- ++(0,-0.2);
        \draw (J) (1.72,0) -- ++(0,0.2) -- ++(0.2,0);   
        \draw (H) (0,1.64) -- ++(0.2,0) -- ++(0,-0.2);
        \coordinate (D) at (1.5,2); 
        \coordinate (E) at (0,2); 
        \begin{scope}
            \clip (0,2) circle (2); 
            \fill[red!30, opacity=0.4, pattern=north west lines]
                (1.92, 1.44) -- 
                (0, 0) -- 
                plot[domain=270:343.74, samples=100]
                    ({0 + 2*cos(\x)}, {2 + 2*sin(\x)})
                -- cycle;
        \end{scope}
        \coordinate (F) at (1.5,0); 
        \draw[blue] (D) circle (2.5); 
        \draw[red] (E) circle (2); 
        \draw[green] (1.5,0) circle (1.5); 
        \coordinate (G) at (1.92,1.44);
        \coordinate (J) at (1.92, 0);
        \coordinate (H) at (0, 1.44);
        \draw[dashed] (G) -- (C);
        \draw[dashed] (E) -- (G);
        \draw[dashed] (1.5,0) -- (G);
        \draw[dotted] (G) -- (J);
        \draw[dotted] (G) -- (H);
        \node at (A) [above left] {A};
        \node at (B) [below right] {B};
        \node at (C) [below left] {C};
        \node at (E) [left] {E};
        \node at (F) [below] {F};
        \node at (G) [above right] {G};
        \node at (H) [left] {H};
        \node at (J) [below] {J};
        \node at (E) [above right] {$b$};
        \node at (F) [above left] {$a$};
        \node at (D) [above right] {$c$};
    \end{tikzpicture}
    \caption{Region D (shaded)}
    \label{fig:fig_d}
\end{figure}
\subsection{\texorpdfstring{Region E: The area of the circular segment subtended by chord $\overline{CG}$ on circle $F$ (Figure \ref{fig:fig_e})}{Region E: The area of the circular segment subtended by chord CG on circle F (Figure 7)}}
We know that $m\angle CFG = 180^{\degree}-2\theta$. Thus, the area of the sector $[\wideparen{CFG}]$ is $\dfrac{\pi a^2}{4} \cdot \dfrac{180^{\degree}-2\theta}{360^{\degree}}$. The area of $[\triangle CFG]$ is $\dfrac{1}{2} GJ \cdot CF$. Substituting values, this gives the area of $[\triangle CFG] = \dfrac{1}{2}(\dfrac{a^2b}{c^2} \cdot \dfrac{a}{2})$, yielding the area below.
\begin{align}
    \text{Area of Region E: Upper-Left Circular Segment of Circle $F$}\nonumber
\end{align}
\begin{align}
    &=\dfrac{\pi a^2}{4} \cdot \dfrac{180^{\degree}-2\theta}{360^{\degree}} -  \dfrac{1}{2}(\dfrac{a^2b}{c^2} \cdot \dfrac{a}{2}) \hypertarget{eq:region e}{}
\end{align}
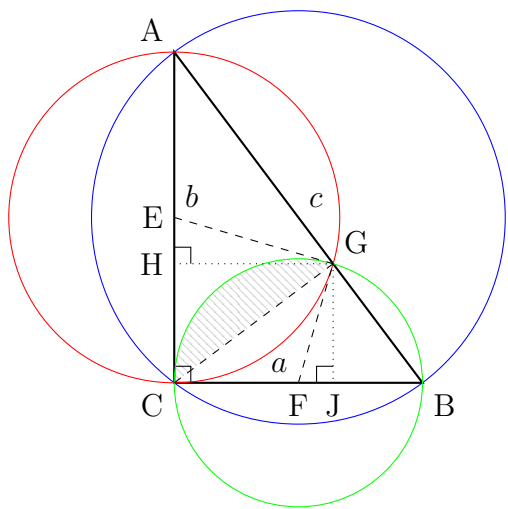
\begin{figure}[H]
    \centering
    \begin{tikzpicture}[scale=1.1]
        \coordinate (C) at (0,0);
        \coordinate (B) at (3,0);
        \coordinate (A) at (0,4);
        \draw[thick] (A) -- (B) -- (C) -- cycle;
        \draw (C) (0,0.2) -- ++(0.2,0) -- ++(0,-0.2);
        \draw (J) (1.72,0) -- ++(0,0.2) -- ++(0.2,0);
        \draw (H) (0,1.64) -- ++(0.2,0) -- ++(0,-0.2);
        \coordinate (D) at (1.5,2); 
        \coordinate (E) at (0,2); 
        \coordinate (F) at (1.5,0); 
        \begin{scope}
            \clip (1.5,0) circle (1.5); 
            \fill[green!30, opacity=0.4, pattern=north west lines]
                (0, 0) -- (1.92, 1.44)
                plot[domain=73.53:180, samples=100]
                    ({1.5 + 1.5*cos(\x)}, {0 + 1.5*sin(\x)})
                -- cycle;
        \end{scope}
        \draw[blue] (D) circle (2.5); 
        \draw[red] (E) circle (2); 
        \draw[green] (1.5,0) circle (1.5); 
        \coordinate (G) at (1.92,1.44);
        \coordinate (J) at (1.92, 0);
        \coordinate (H) at (0, 1.44);
        \draw[dashed] (G) -- (C);
        \draw[dashed] (E) -- (G);
        \draw[dashed] (1.5,0) -- (G);
        \draw[dotted] (G) -- (J);
        \draw[dotted] (G) -- (H);
        \node at (A) [above left] {A};
        \node at (B) [below right] {B};
        \node at (C) [below left] {C};
        \node at (E) [left] {E};
        \node at (F) [below] {F};
        \node at (G) [above right] {G};
        \node at (H) [left] {H};
        \node at (J) [below] {J};        
        \node at (E) [above right] {$b$};
        \node at (F) [above left] {$a$};
        \node at (D) [above right] {$c$};      
    \end{tikzpicture}
    \caption{Region E (shaded)}
    \label{fig:fig_e}
\end{figure}

\subsection{\texorpdfstring{Region F: The area of the circular segment $\wideparen{BG}$ subtended by chord $\overline{BG}$ on circle $F$ (Figure \ref{fig:fig_f})}{Region F: The area of the circular segment subtended by chord BG on circle F (Figure 8)}}
We know that $m\angle BFG = 2\theta$. Thus, the area of the sector $[\wideparen{BFG}]$ for Circle $F$ is $\dfrac{\pi a^2}{4} \cdot \dfrac{2\theta}{360^{\degree}}$. The area of the  $[\triangle BFG]$ is $\dfrac{1}{2} GJ \cdot FB$. Substituting values, we see that the area of $[\triangle BFG] = \dfrac{1}{2}(\dfrac{a^2b}{c^2} \cdot \dfrac{a}{2})$. Thus, we get the area below.
\begin{align}
    \text{Area of Region F: Upper Right Circular Segment of Circle $F$}\nonumber
\end{align}
\begin{align}
    &= \dfrac{\pi a^2}{4} \cdot \dfrac{2\theta}{360^{\degree}} - \dfrac{1}{2}(\dfrac{a^2b}{c^2} \cdot \dfrac{a}{2}) \hypertarget{eq:region f}{}
\end{align}
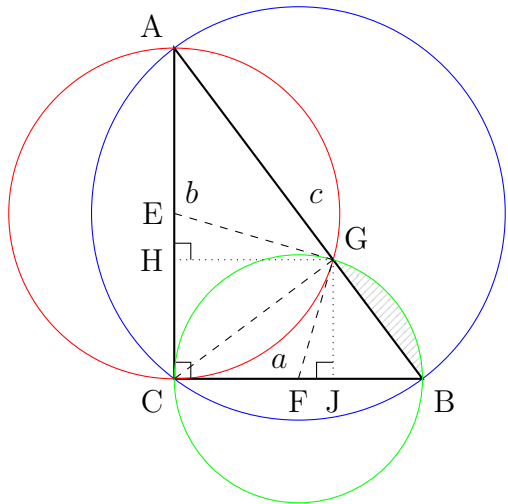
\begin{figure}[H]
    \centering
    \begin{tikzpicture}[scale=1.1]
        \coordinate (C) at (0,0);
        \coordinate (B) at (3,0);
        \coordinate (A) at (0,4);
        \draw[thick] (A) -- (B) -- (C) -- cycle;
        \draw (C) (0,0.2) -- ++(0.2,0) -- ++(0,-0.2);
        \draw (J) (1.72,0) -- ++(0,0.2) -- ++(0.2,0);
        \draw (H) (0,1.64) -- ++(0.2,0) -- ++(0,-0.2);
        \coordinate (D) at (1.5,2); 
        \coordinate (E) at (0,2); 
        \coordinate (F) at (1.5,0); 
        \begin{scope}
            \clip (1.5,0) circle (1.5); 
            \fill[green!30, opacity=0.4, pattern=north east lines]
                (3, 0) -- (1.92, 1.44)
                plot[domain=73.53:0, samples=100]
                    ({1.5 + 1.5*cos(\x)}, {0 + 1.5*sin(\x)})
                -- cycle;
        \end{scope}
        \draw[blue] (D) circle (2.5); 
        \draw[red] (E) circle (2); 
        \draw[green] (1.5,0) circle (1.5); 
        \coordinate (G) at (1.92,1.44);
        \coordinate (J) at (1.92, 0);
        \coordinate (H) at (0, 1.44);
        \draw[dashed] (G) -- (C);
        \draw[dashed] (E) -- (G);
        \draw[dashed] (1.5,0) -- (G);
        \draw[dotted] (G) -- (J);
        \draw[dotted] (G) -- (H);
        \node at (A) [above left] {A};
        \node at (B) [below right] {B};
        \node at (C) [below left] {C};
        \node at (E) [left] {E};
        \node at (F) [below] {F};
        \node at (G) [above right] {G};
        \node at (H) [left] {H};
        \node at (J) [below] {J};
        \node at (E) [above right] {$b$};
        \node at (F) [above left] {$a$};
        \node at (D) [above right] {$c$};
    \end{tikzpicture}
    \caption{Region F (shaded)}
    \label{fig:fig_f}
\end{figure}

\section{\texorpdfstring{To Find The Area of A Semicircle}{The Area of the Semicircle}} \label{sec:main body}
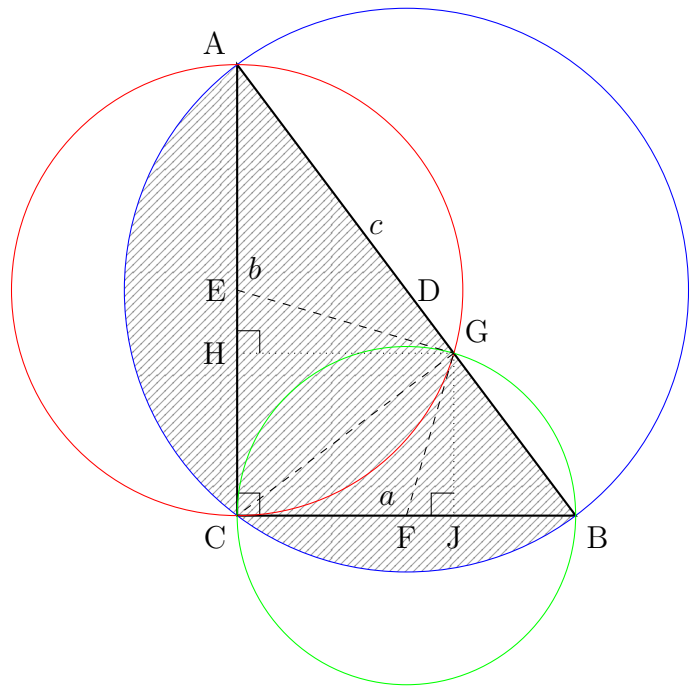
\begin{figure}[H]
    \centering
    \begin{tikzpicture}[scale=1.5]
        \coordinate (C) at (0,0); 
        \coordinate (B) at (3,0);
        \coordinate (A) at (0,4);
        \draw[thick] (A) -- (B) -- (C) -- cycle;
        \draw (C) (0,0.2) -- ++(0.2,0) -- ++(0,-0.2);
        \draw (J) (1.72,0) -- ++(0,0.2) -- ++(0.2,0);
        \draw (H) (0,1.64) -- ++(0.2,0) -- ++(0,-0.2);
        \coordinate (D) at (1.5,2); 
        \begin{scope}
            \clip (1.5,2) circle (2.5); 
            \fill[yellow!30, opacity=0.6, pattern=north east lines]
                (A) -- (B)
                arc (323.130102354:143.130102354:2.5)
                -- cycle;
        \end{scope}
        \coordinate (E) at (0,2); 
        \coordinate (F) at (1.5,0); 
        \draw[blue] (D) circle (2.5); 
        \draw[red] (E) circle (2); 
        \draw[green] (1.5,0) circle (1.5); 
        \coordinate (G) at (1.92,1.44);
        \coordinate (J) at (1.92, 0);
        \coordinate (H) at (0, 1.44);
        \draw[dashed] (G) -- (C);
        \draw[dashed] (E) -- (G);
        \draw[dashed] (1.5,0) -- (G);
        \draw[dotted] (G) -- (J);
        \draw[dotted] (G) -- (H);
        \node at (A) [above left] {A};
        \node at (B) [below right] {B};
        \node at (C) [below left] {C};
        \node at (D) [right] {D};
        \node at (E) [left] {E};
        \node at (F) [below] {F};
        \node at (G) [above right] {G};
        \node at (H) [left] {H};
        \node at (J) [below] {J};
        \node at (E) [above right] {$b$};
        \node at (F) [above left] {$a$};
        \node at (1.23, 2.4) [above] {$c$};
    \end{tikzpicture}
    \caption{Semicircle of Circle $D$ containing $\triangle ABC$ (shaded)}
    \label{fig:Semicircle}
\end{figure}
    The area of the semicircle over the hypotenuse, as depicted in Figure \ref{fig:Semicircle}, can also be expressed as the sum of the areas of the semicircles over the triangle’s legs (circles E and F), together with the circular segments from circle C labeled as Regions A and B, subtracting the circular segments from circles E and F (namely, Regions C through F).

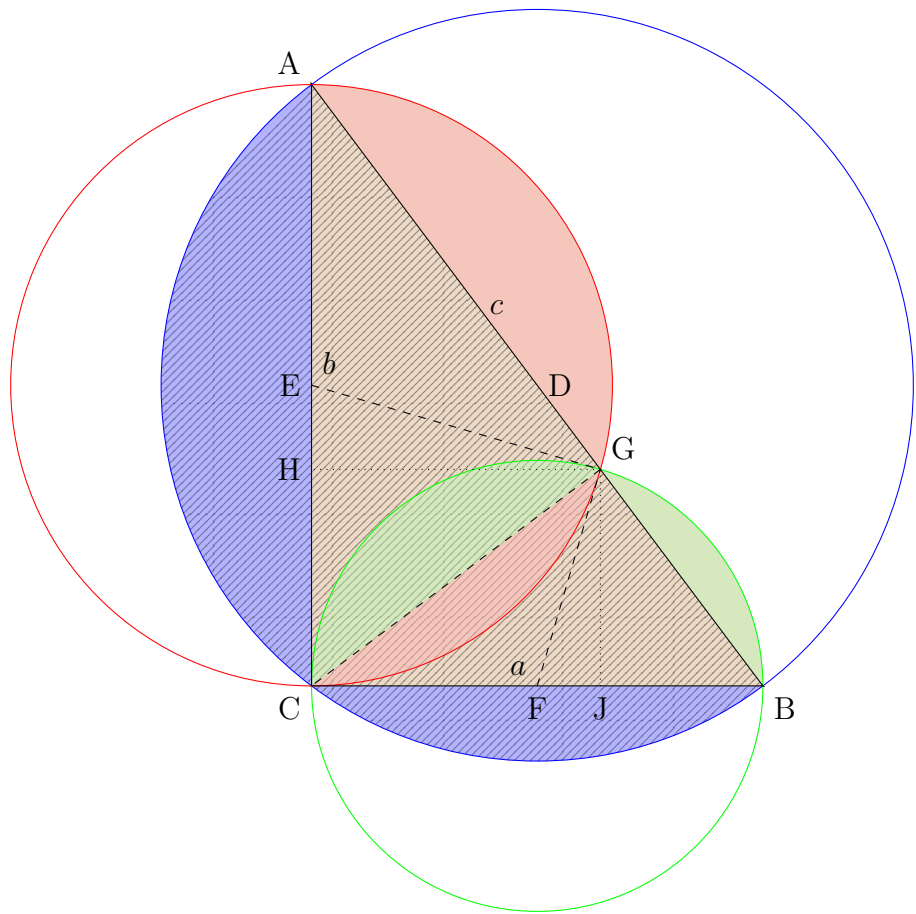
\begin{figure}[H]
    \centering
    \begin{tikzpicture}[scale=2]
        \coordinate (C) at (0,0);
        \coordinate (B) at (3,0);
        \coordinate (A) at (0,4);
        \draw[thick] (A) -- (B) -- (C) -- cycle;
        \draw (C) (0,0.2) -- ++(0.2,0) -- ++(0,-0.2);
        \draw (J) (1.72,0) -- ++(0,0.2) -- ++(0.2,0); 
        \draw (H) (0,1.64) -- ++(0.2,0) -- ++(0,-0.2);
        \coordinate (D) at (1.5, 2);
        \begin{scope}
            \clip (D) circle (2.5); 
            \fill[blue!30, opacity=1] (C) -- (B) arc (306.87:233.13:2.5) -- cycle;
        \end{scope}
        \begin{scope}
            \clip (D) circle (2.5); 
            \fill[blue!30, opacity=1] (A) -- (C) arc (233.13:126.87:2.5) -- cycle;
        \end{scope}
        \begin{scope}
            \clip (0,2) circle (2);
            \fill[brown!30, opacity=1]
                (0,4) -- (0,0)
                plot[domain=-90:90, samples=100]
                    ({0 + 2*cos(\x)}, {2 + 2*sin(\x)})
                -- cycle;
        \end{scope}
        \begin{scope}
            \clip (1.5,0) circle (1.5); 
            \fill[green!30, opacity=1]
                (0, 0) -- (1.92, 1.44)
                plot[domain=73.53:180, samples=100]
                    ({1.5 + 1.5*cos(\x)}, {0 + 1.5*sin(\x)})
                -- cycle;
                \clip (1.5,0) circle (1.5); 
        \end{scope}

        \begin{scope}
            \clip (1.5,0) circle (1.5); 
            \fill[brown!30, opacity=1]
                (0, 0) -- (1.92, 1.44)
                plot[domain=0:180, samples=100]
                    ({1.5 + 1.5*cos(\x)}, {0 + 1.5*sin(\x)})
                -- cycle;
        \end{scope}
        \begin{scope}
            \clip (1.5,2) circle (2.5); 
            \fill[yellow!30, opacity=0.6, pattern=north east lines]
                (A) -- (B)
                arc (323.130102354:143.130102354:2.5)
                -- cycle;
        \end{scope}
        \begin{scope}
            \clip (0,2) circle (2); 
            \fill[red!30, opacity=0.5]
                (1.92, 1.44) -- 
                (0, 0) -- 
                plot[domain=270:343.74, samples=100]
                    ({0 + 2*cos(\x)}, {2 + 2*sin(\x)})
                -- cycle;
        \end{scope}
        \begin{scope}
            \clip (1.5,0) circle (1.5); 
            \fill[green!30, opacity=0.4]
                (0, 0) -- (1.92, 1.44)
                plot[domain=73.53:180, samples=100]
                    ({1.5 + 1.5*cos(\x)}, {0 + 1.5*sin(\x)})
                -- cycle;

                \clip (1.5,0) circle (1.5); 
        \end{scope}

        \coordinate (E) at (0,2); 
        
        \begin{scope}
            \clip (0,2) circle (2);
            \fill[red!30, opacity=0.5]
                (0,4) -- (1.92,1.44)
                plot[domain=-16.26:90, samples=100]
                    ({0 + 2*cos(\x)}, {2 + 2*sin(\x)})
                -- cycle;
        \end{scope}

        \coordinate (F) at (1.5,0); 
        \begin{scope}
            \clip (1.5,0) circle (1.5); 
            \fill[green!30, opacity=0.4]
                (3, 0) -- (1.92, 1.44)
                plot[domain=73.53:0, samples=100]
                    ({1.5 + 1.5*cos(\x)}, {0 + 1.5*sin(\x)})
                    
                -- cycle;
            pattern
        \end{scope}

        \draw[blue] (D) circle (2.5); 
        \draw[red] (E) circle (2); 
        \draw[green] (1.5,0) circle (1.5); 
        
        \coordinate (G) at (1.92,1.44);
        \coordinate (J) at (1.92, 0);
        \coordinate (H) at (0, 1.44);
        \draw[dashed] (G) -- (C);
        \draw[dashed] (E) -- (G);
        \draw[dashed] (1.5,0) -- (G);
        \draw[dotted] (G) -- (J);
        \draw[dotted] (G) -- (H);
        \draw (A) -- (B) -- (C) -- cycle;
        \node at (A) [above left] {A};
        \node at (B) [below right] {B};
        \node at (C) [below left] {C};
        \node at (D) [right] {D};
        \node at (E) [left] {E};
        \node at (F) [below] {F};
        \node at (G) [above right] {G};
        \node at (H) [left] {H};
        \node at (J) [below] {J};
        
        \node at (E) [above right] {$b$};
        \node at (F) [above left] {$a$};
        \node at (1.23, 2.4) [above] {$c$};
        
    \end{tikzpicture}
    \caption{Conglomerate of Regions A, B, C, D, E, F, and Semicircles of Circles $E$ and $F$}
    \label{fig:Semicircle_Conglomerate}
\end{figure}

    Figure \ref{fig:Semicircle_Conglomerate} illustrates this decomposition: the brown regions correspond to the leg semicircles, while the overlaid red and green sectors represent the circular segments that are subtracted. The resulting shaded area reflects the net area obtained after performing the indicated additions and subtractions.

\subsection{Formal Area Decomposition}

We now consolidate the area expressions derived in Section 5 into a rigorous decomposition identity. This identity asserts that the area of the semicircle constructed on the hypotenuse $AB$ is fully reconstructible via the semicircles on the legs and the net contribution of specific circular segments.

Let Regions A–F denote the circular segments defined in Figures \ref{fig:fig_a}–\ref{fig:fig_f}, and let $S_c$, $S_a$, and $S_b$ denote the semicircles constructed over the hypotenuse and legs $a$ and $b$. Then:
\begin{align}
\text{
$
\boxed{[S_c] = [S_a] + [S_b] + [\hyperlink{eq:region a}{A}] + [\hyperlink{eq:region b}{B}]  - [\hyperlink{eq:region c}{C}] - [\hyperlink{eq:region d}{D}] - [\hyperlink{eq:region e}{E}] - [\hyperlink{eq:region f}{F}]}
$
}
\end{align}

\begin{proof}
As established in the construction of Section 3 and detailed in Section 5, the semicircle on the hypotenuse $ AB $ overlaps in part with the semicircles constructed over legs $ AC $ and $ BC $. The regions labeled A and B are exterior circular segments formed between the hypotenuse arc and the legs, while Regions C through F represent areas of overlap or excess arising from those leg-based semicircles.

To obtain the area of the hypotenuse semicircle from the component parts:
\begin{itemize}
    \item Begin with the semicircles on legs $ AC $ and $ BC $
    \item Add the two external circular segments labeled Region $\hyperlink{eq:region a}{A}$ (over arc $\wideparen{AC}$) and Region $\hyperlink{eq:region b}{B}$ (over arc $\wideparen{BC}$) that lie within the hypotenuse’s arc but outside the leg arcs.
    \item  Subtract Regions $\hyperlink{eq:region c}{C}$ through $\hyperlink{eq:region f}{F}$, which represent overlapping or extraneous areas: intersections between the hypotenuse semicircle and the leg-based semicircles.
\end{itemize}

\vspace{1em}
\noindent\textbf{Formality of Region Decomposition.}
To elevate the decomposition beyond visual intuition, we now rigorously justify that each region included in the area equation is both well-defined and non-overlapping.

\begin{lemma}[Region A Congruence]
Let $\wideparen{AC}$ denote the arc of Circle $D$ spanning segment $AC$, and let $\wideparen{AG}$ and $\wideparen{CG}$ denote arcs in Circle $E$. Then, Region A (as labeled in Figure~\ref{fig:fig_a}) is equal in area to the union of the two circular segments bounded by $\wideparen{AG}$ and $\wideparen{CG}$ minus the area of $\triangle AGC$. This congruence follows from rotational symmetry and the equality of corresponding arcs and radii in the construction.
\end{lemma}

\begin{lemma}[Region Partition Exhaustiveness]
The regions labeled A–F, along with the semicircles over legs $AC$ and $BC$, are disjoint and exhaustive with respect to the semicircle constructed on the hypotenuse $AB$. Each region lies entirely within Circle $D$ and is bounded by distinct arcs and chords. Their union, with appropriate signed inclusion or subtraction, exactly tiles the hypotenuse semicircle without overlap.
\end{lemma}

\begin{remark}
This decomposition relies solely on the construction of circles over triangle sides and angle-based partitions; it does not presuppose the Pythagorean Theorem. Trigonometric identities (e.g., $\sin(2\theta)$) and geometric symmetries underlie the area relationships, ensuring that the decomposition is internally self-consistent and grounded in first principles. See footnotes for references. 
\end{remark}
 
 This establishes a purely geometric identity that is both disjoint and exhaustive grounded in spatial congruence.
\end{proof}

\subsection{\texorpdfstring{Derivation of the Area [$S_c$]}{Derivation of the Area of the Semicircle of Circle D}}
To complete the argument, we will express the total area of the semicircle [$S_c$] over side $AB$ in two distinct ways: once directly through $\pi r^2$, and once as a combination of segments and semicircles from the side circles. Equating the two expressions yields the Pythagorean identity.

The area of a semicircle of circle $D$ (Figure \ref{fig:Semicircle}), $\boxed{S_c = \dfrac{\pi c^2}{8}}$ when calculated using the formula $\pi r^2$, as shown in Figure \ref{fig:Semicircle}.

\begin{align}
        \text{\boxed{S_c = 
        \frac{\pi a^2}{8} +
        \frac{\pi b^2}{8} + 
        [\hyperlink{eq:region a}{A}] + [\hyperlink{eq:region b}{B}] - [\hyperlink{eq:region c}{C}] - [\hyperlink{eq:region d}{D}] - [\hyperlink{eq:region e}{E}] - [\hyperlink{eq:region f}{F}]
        }}\nonumber
\end{align}
\begin{align}
    \Downarrow\nonumber
\end{align}
\begin{align} =
        (\dfrac{\pi c^2}{4} \cdot \dfrac{180^{\degree} - 2\theta}{360^{\degree}} - \dfrac{ab}{4})\; + \;(\dfrac{\pi c^2}{4} \cdot \dfrac{2\theta}{360^{\degree}} - \dfrac{ab}{4}) \\ 
        -\:(\dfrac{\pi b^2}{4} \cdot \dfrac{180^{\degree} - 2\theta}{360^{\degree}} - \dfrac{1}{2}(\dfrac{ab^2}{c^2} \cdot \dfrac{b}{2}))\; - \;(\dfrac{\pi b^2}{4} \cdot \dfrac{2\theta}{360^{\degree}} - \dfrac{1}{2}(\dfrac{ab^2}{c^2} \cdot \dfrac{b}{2}))  \nonumber \\ 
        -\:(\dfrac{\pi a^2}{4} \cdot \dfrac{180^{\degree}-2\theta}{360^{\degree}} - \dfrac{1}{2}(\dfrac{a^2b}{c^2} \cdot \dfrac{a}{2}))\; - \;(\dfrac{\pi a^2}{4} \cdot \dfrac{2\theta}{360^{\degree}} - \dfrac{1}{2}(\dfrac{a^2b}{c^2} \cdot \dfrac{a}{2})) \nonumber \\ 
        +\:\dfrac{\pi a^2}{8}  +\: \dfrac{\pi b^2}{8} \nonumber
\end{align}
Reorganizing terms, we arrive at:
\begin{align}
    =\: \dfrac{\pi c^2}{4} \cdot (\dfrac{180^{\degree} - 2\theta}{360^{\degree}} + \dfrac{2\theta}{360^{\degree}}) - \dfrac{ab}{2} \\ 
        -\:\dfrac{\pi b^2}{4} \cdot (\dfrac{180^{\degree} - 2\theta}{360^{\degree}} +  \dfrac{2\theta}{360^{\degree}}) + (\dfrac{1}{2}(\dfrac{ab^2}{c^2}\cdot\dfrac{b}{2})\: + \: \dfrac{1}{2}(\dfrac{ab^2}{c^2}\cdot\dfrac{b}{2}))\; \nonumber \\ 
        -\:\dfrac{\pi a^2}{4} \cdot (\dfrac{180^{\degree}-2\theta}{360^{\degree}} + \dfrac{2\theta}{360^{\degree}}) + (\dfrac{1}{2}(\dfrac{a^2b}{c^2}\cdot\dfrac{a}{2})\: + \: \dfrac{1}{2}(\dfrac{a^2b}{c^2}\cdot\dfrac{a}{2})) \nonumber \\ 
        +\:\dfrac{\pi a^2}{8}  +\: \dfrac{\pi b^2}{8} \nonumber
\end{align}
We can now change $\dfrac{180 - 2\theta + 2\theta}{360}$ into $\dfrac{1}{2}$, removing the variable $\theta$, and simply add together the $\dfrac{1}{2}(\cdot\cdot\cdot)$ terms. Although the angle $\theta$ appears explicitly in several intermediate area expressions, all such dependencies cancel exactly in the final decomposition, confirming that the proof does not rely on specific angle measures and remains valid for all right triangles:
\begin{align}
    =\: \dfrac{\pi c^2}{4} \cdot \dfrac{1}{2} - \dfrac{ab}{2} 
        -\:\dfrac{\pi b^2}{4} \cdot \dfrac{1}{2} + \dfrac{ab^2}{c^2} \cdot \dfrac{b}{2}\;
        -\:\dfrac{\pi a^2}{4} \cdot \dfrac{1}{2} + \dfrac{a^2b}{c^2} \cdot \dfrac{a}{2} \\ 
        +\:\dfrac{\pi a^2}{8}  +\: \dfrac{\pi b^2}{8} \nonumber
\end{align}
\begin{align}
    =\: \dfrac{\pi c^2}{8} - \dfrac{ab}{2}
        -\:\dfrac{\pi b^2}{8}  + \dfrac{ab^3}{2c^2}\; 
        -\:\dfrac{\pi a^2}{8} + \dfrac{a^3b}{2c^2}
        +\:\dfrac{\pi a^2}{8}  +\: \dfrac{\pi b^2}{8} 
\end{align}
Now, the terms begin to cancel out neatly:
\begin{align}
    &=\:\dfrac{\pi c^2}{8}\: - \:\dfrac{ab}{2}\: + \:\dfrac{ab^3}{2c^2}\: + \:\dfrac{a^3b}{2c^2} \\ \nonumber \\
    &=\dfrac{\pi c^2}{8} - \dfrac{ab}{2} + \dfrac{a^3b}{2c^2} + \dfrac{ab^3}{2c^2} 
\end{align}
Now we factor out $ab$ from $a^3b + ab^3$:
\begin{align}
    \text{$
    = \dfrac{\pi c^2}{8} - \dfrac{ab}{2} + \dfrac{ab(a^2+b^2)}{2c^2}
    $}
\end{align}
After simplification, we see the area of the semicircle of Circle D that contains $\triangle ABC$, $\: \boxed{S_c = \dfrac{\pi c^2}{8} - \dfrac{ab}{2} + \dfrac{ab(a^2+b^2)}{2c^2}}$.
\subsection{The Final Stretch}
We now set these two areas equal to each other, as they represent the same area and are thus equal:
\begin{align}
    S_c = S_c \\ \nonumber \\
    \dfrac{\pi c^2}{8} - \dfrac{ab}{2} + \dfrac{ab(a^2+b^2)}{2c^2} &= \dfrac{\pi c^2}{8} \\  \nonumber \\
    - \dfrac{ab}{2} + \dfrac{ab(a^2+b^2)}{2c^2} &= 0 \\  \nonumber \\
    \dfrac{ab(a^2+b^2)}{2c^2} &= \dfrac{ab}{2}
\end{align}

Multiplying by ${2c^2}$:
\begin{align}
    ab(a^2+b^2) &= c^2ab
\end{align}
Divide both sides by $ab$:
\begin{align}
    \dfrac{ab(a^2+b^2)}{ab} &= \frac{c^2ab}{ab} \\ \nonumber \\ 
    a^2 + b^2 &= c^2
\end{align}
Thus, we finally prove $\boxed{a^2 + b^2 = c^2}$.

\section{Conclusion}

This proof of the Pythagorean theorem provides a novel approach using circles to circumscribe the sides of a right triangle. By analyzing the areas of the segments formed by these circles and their intersections, we derive the classic relationship $a^2 + b^2 = c^2$. The proof combines geometric construction with trigonometric identities and area comparisons in a novel configuration not previously documented in geometric literature. In our manual search of The Pythagorean Proposition by Loomis \cite{loomis}, which catalogs hundreds of known proofs and manual review of the extensive online repository maintained by Bogomolny \cite{bogomolny}, no proof was found to use this particular decomposition strategy involving overlapping circular segments, circles, and basic trigonometry. This suggests that the method outlined here may constitute a novel contribution to the known corpus of geometric and trigonometric proofs of the Pythagorean Theorem. 

This proof reveals that the square relationship emerges naturally from the structure of angle-determined sectors, tying together angular geometry and area in a way not visible in commonly found square-based or rearrangement proofs. While algebraic in execution, this proof’s structure is inherently geometric, grounded in symmetry and spatial decomposition rather than coordinate or vector methods. 

Further exploration may expand upon this diagram and challenge younger generations to find innovative proofs within new and old diagrams that will shape the mathematical world. This paper hopefully serves as an inspiration for the next generation to seek novelty in what seems exhausted, search for what is hidden in plain sight, build upon the works of others that came before, and continue to challenge oneself in new and exciting ways.

\section{Related Work}
\subsection{Statement of Novelty}
A manual review of 370+ entries in Loomis and the full Cut-the-Knot database revealed no prior use of this geometric configuration or decomposition. This suggests that the proof presented here is new to the documented mathematical literature.

\subsection{Another Method} \label{sec:another method}

In a previous footnote, we mentioned how the proof can omit the use of trigonometry. In fact, we can remove the use of $sin(2\theta)$ and any use or understanding of $sin(\alpha)$ or $cos(\alpha)$. 

In this proof, we utilize trigonometry to calculate the lengths of altitudes $GH$ and $GJ$. We then use these altitudes to calculate the areas of 4 triangles: $\triangle CEG$, $\triangle AEG$, $\triangle CFG$,  $\triangle GFB$. 

Instead of using trigonometry to calculate these areas, we can use similar triangles. 
\begin{figure}[H]
    \centering
    \begin{tikzpicture}[scale=1.56]
        \coordinate (C) at (0,0);
        \coordinate (B) at (3,0);
        \coordinate (A) at (0,4);
        \coordinate (J) at (1.92, 0);
        \coordinate (H) at (0, 1.44);
        \draw[thick] (A) -- (B) -- (C) -- cycle;
        \draw (C) (0,0.2) -- ++(0.2,0) -- ++(0,-0.2);
        \draw (J) (1.72,0) -- ++(0,0.2) -- ++(0.2,0);
        \draw (H) (0,1.64) -- ++(0.2,0) -- ++(0,-0.2);
        \coordinate (D) at (1.5,2); 
        \coordinate (E) at (0,2); 
        \coordinate (F) at (1.5,0); 
        \draw[blue] (D) circle (2.5); 
        \draw[red] (E) circle (2); 
        \draw[green] (1.5,0) circle (1.5); 
        \coordinate (G) at (1.92,1.44);
        
        \draw[dashed] (G) -- (C);
        \draw[dashed] (E) -- (G);
        \draw[dashed] (1.5,0) -- (G);
        \draw[dotted] (G) -- (J);
        \draw[dotted] (G) -- (H);
        \draw[dashed] (C) -- (D);
        \node at (A) [above left] {A};
        \node at (B) [below right] {B};
        \node at (C) [below left] {C};
        \node at (D) [right] {D};
        \node at (E) [left] {E};
        \node at (F) [below] {F};
        \node at (G) [above right] {G};
        \node at (H) [left] {H};
        \node at (J) [below] {J};
        \node at (E) [above right] {$b$};
        \node at (F) [above left] {$a$};
        \node at (1.1, 2.7) [right] {$c$};
    \end{tikzpicture}
    \caption{Circumscribed Construction including Line Segments $\overline{GH}$ and $\overline{GJ}$ For Alternate Proof}
    \label{fig:another_method_diagram}
\end{figure}
We have already established that $\angle CGA$ and $\angle CGB$ are right in \nameref{sec:Construction} and Figure \ref{fig:init_diagram}. Therefore, we observe that $\triangle AGC \sim \triangle ACB \sim \triangle CGB$ via Angle-Angle similarity. Since these triangles are similar, we can find the lengths of their sides. Since $AC = b$, we determine that $AG = \dfrac{b}{c}\cdot b$, $CG = \dfrac{b}{c} \cdot a = \dfrac{a}{c}\cdot b$, and $BG = \dfrac{a}{c} \cdot a$. From this, we can find the area of the 4 triangles from earlier without using trignometry. Note that $AE = CE$ and $CF = FB$ because all radii of a circle are equal: thus, since all triangles with the same base and same height have equal area, we determine that $[\triangle AEG] = \dfrac{1}{2}\cdot [\triangle ACG], \: [\triangle CEG] = \dfrac{1}{2}\cdot [\triangle ACG],  \: [\triangle CFG] = \dfrac{1}{2}\cdot [\triangle CGB],\: [\triangle BGF] = \dfrac{1}{2}\cdot [\triangle CGB]$

\begin{align} 
    \text{$[\triangle CEG] = \dfrac{1}{2}\cdot[ \triangle CGA] = \dfrac{1}{2}\cdot( \dfrac{1}{2} \cdot AG \cdot CG) = \dfrac{1}{2}\cdot( \dfrac{1}{2} \cdot \dfrac{b^2}{c} \cdot \dfrac{ab}{c}) = \dfrac{1}{4}\cdot\dfrac{ab^3}{c^2}$}
    \\
    \text{$[\triangle AEG] = \dfrac{1}{2}\cdot[ \triangle CGA] = \dfrac{1}{2}\cdot( \dfrac{1}{2} \cdot AG \cdot CG) = \dfrac{1}{2}\cdot( \dfrac{1}{2} \cdot \dfrac{b^2}{c} \cdot \dfrac{ab}{c}) = \dfrac{1}{4}\cdot\dfrac{ab^3}{c^2}$}
    \\
    \text{$[\triangle CFG] = \dfrac{1}{2}\cdot[ \triangle CGB] = \dfrac{1}{2}\cdot( \dfrac{1}{2} \cdot CG \cdot GB) = \dfrac{1}{2}\cdot( \dfrac{1}{2} \cdot \dfrac{ab}{c} \cdot \dfrac{a^2}{c}) = \dfrac{1}{4}\cdot\dfrac{a^3b}{c^2}$}
    \\
    \text{$[\triangle GFB] = \dfrac{1}{2}\cdot[ \triangle CGB] = \dfrac{1}{2}\cdot( \dfrac{1}{2} \cdot CG \cdot GB) = \dfrac{1}{2}\cdot( \dfrac{1}{2} \cdot \dfrac{ab}{c} \cdot \dfrac{a^2}{c}) = \dfrac{1}{4}\cdot\dfrac{a^3b}{c^2}$}
\end{align}
Thus, without using any trigonometric methods, we are able to find these areas and can complete the proof using the same procedure and algebra in Section \ref{sec:main body}.

\section{Acknowledgments}

I would like to express gratitude to Lily Marin Parsapour, Leonidas Boukas, and Major Lockhart Matthews for their help in reviewing and verifying my proof. I also acknowledge the assistance of Claude 3.7, an AI assistant developed by Anthropic, and ChatGPT-4o, an AI assistant developed by OpenAI, who both provided support in drafting the template and structure for this paper. While all mathematical work, proofs, and analysis are my own original contributions, Claude AI and ChatGPT helped with organizing the content and formatting the document.

\end{document}